\documentclass[12pt]{article}
\usepackage{amssymb}
\begin{document}

\def\kbar{{\mathchar'26\mkern-9muk}}  
\def\bra#1{\langle #1 \vert}
\def\ket#1{\vert #1 \rangle}
\def\vev#1{\langle #1 \rangle}
\def\tr{\mbox{Tr}\,}
\def\ad{\mbox{ad}\,}
\def\ker{\mbox{Ker}\,}
\def\im{\mbox{Im}\,}
\def\der{\mbox{Der}\,}
\def\ad{\mbox{ad}\,}
\def\b#1{{\mathbb #1}}
\def\op{{\scriptstyle\rm op}}
\def\e{{\scriptstyle\rm e}}
\def\c#1{{\cal #1}}
\def\pt{\partial_t}
\def\px{\partial_1}
\def\bpx{\bar\partial_1}
\def\pr{\partial_{R1}}
\newcommand{\mat}[4]{\left(\!\begin{array}{cc}#1&#2\\#3&#4\end{array}\!\right)}
\def\m@th{\mathsurround=0pt}
\def\eqalign#1{\null\,\vcenter{\openup 3pt \m@th
\ialign{\strut\hfil$\displaystyle{##}$&$\displaystyle{{}##}$\hfil
\crcr#1\crcr}}\,}
\newcommand{\sect}[1]{\setcounter{equation}{0}\section{#1}}
\renewcommand{\theequation}{\thesection.\arabic{equation}}
\newcommand{\subsect}[1]{\setcounter{equation}{0}\subsection{#1}}
\renewcommand{\theequation}{\thesection.\arabic{equation}}
\newcommand{\be}{\begin{equation}}
\newcommand{\ee}{\end{equation}}
\newcommand{\bea}{\begin{eqnarray}}
\newcommand{\eea}{\end{eqnarray}}

\title{The Geometry of a $q$-Deformed Phase Space}

\author{B.L. Cerchiai$\strut^{1,2}$ \,
        R. Hinterding$\strut^{1,2}$ \\ J. Madore$\strut^{2,3}$, 
        and \ J. Wess$\strut ^{1,2}$
        \and
        $\strut^1$Sektion Physik, Ludwig-Maximilian Universit\"at,\\
        Theresienstra\ss e 37, D-80333 M\"unchen
        \and
        $\strut^2$Max-Planck-Institut f\"ur Physik\\
        (Werner-Heisenberg-Institut)\\
        F\"ohringer Ring 6, D-80805 M\"unchen
        \and
        $\strut^3$Laboratoire de Physique Th\'eorique et Hautes Energies\\
        Universit\'e de Paris-Sud, B\^atiment 211, F-91405 Orsay
        }

\date{}

\maketitle

\abstract{The geometry of the $q$-deformed line is studied. A real
differential calculus is introduced and the associated algebra of
forms represented on a Hilbert space. It is found that there is a
natural metric with an associated linear connection which is of zero
curvature. The metric, which is formally defined in terms of
differential forms, is in this simple case identifiable as an observable.}

\vfill
\noindent
LMU Preprint  98/08
\medskip
\eject

\parskip 4pt plus2pt minus2pt

\sect{Introduction and motivation}

There is a particularly simple noncommutative geometry, the
1-dimensional $q$-deformed euclidean space~\cite{FadResTak89,
CarSchWat91, FicLorWes96} $\b{R}_q^1$ which can be completely analyzed
from almost every point of view. Although this `space' has but one
`dimension' and therefore there are no curvature effects, the
corresponding algebra is most conveniently generated using elements
which in the commutative limit correspond to coordinates in which the
metric does not take its canonical flat form. If one writes, for
example, the line element which describes the distance along the
$y$-axis $ds^2 = dy^2$ using the coordinate $x = e^y$ then one must
write $ds^2 = x^{-2} dx^2$; the metric has component $g_{11} = x^{-2}$. 
It has been argued previously~\cite{DimMad96} that to within
a scale factor there is essentially a unique metric consistent with
the noncommutative structure of an algebra. We shall see this clearly
is the present example.  We shall give a description of this metric in
all detail since it is one of the rare cases in which the general
formalism can be understood in terms of simple physical
observables. In this section we shall give a brief review of the
description of the differential structure of a noncommutative `space'
from the point of view of differential forms and from the `dual' point
of view of twisted derivations. In Section~2, after a few introductory
remarks concerning the algebras $\b{C}^n_q$ and $\b{R}^n_q$ for
general $n$, we describe the algebra $\b{R}^1_q$. In Section~3 we
introduce two conjugate differential calculi over this algebra and in
Section~4 we propose a construction of a real differential
calculus. In Section~5 we discuss a `dual' point of view using twisted
derivations. In Section~6 we briefly mention integration. In Section~7
we discuss the geometry of $\b{R}^1_q$ using the unique local
metric. In Section~8 we introduce Yang-Mills fields and in Section~9
we discuss the Schr\"odinger equation and the Klein-Gordon equation.
In Section~10 we define an associated phase space~\cite{FicLorWes96}
and briefly discuss the harmonic oscillator.  The final section is
devoted to a discussion of the effects of choosing an alternative
non-local metric. Implicitly this metric has been used 
before~\cite{LorRufWes97}.

Noncommutative geometry is geometry which is described by an
associative algebra $\c{A}$ which is usually but not essentially
noncommutative and in which the set of points, if it exists at all, is
relegated to a secondary role. For a thorough exposition of the
subject we refer to the book by Connes~\cite{Con94}; for a gentle
introduction we refer to Madore~\cite{Mad95} or to
Landi~\cite{Lan97}. We shall be exclusively interested here in
algebras which are in some sense deformations of algebras of smooth
functions over a manifold.  A differential calculus over $\c{A}$ is
another associative algebra $\Omega^*(\c{A})$, with a differential
$d$, which plays the role of the de~Rham differential calculus and
must tend to this calculus in the commutative limit. The differential
calculus is what gives structure to the set of `points'. It determines
the `dimension' for example. It would determine the number of nearest
neighbours in the case of a lattice. Over a given $\c{A}$ one can
construct many differential calculi and the one which one choses
depends evidently on the limit manifold one has in mind. There are
many ways one can construct differential calculi. Historically the
first construction~\cite{Con94} was based on an operator which played
in some sense the role of the Dirac operator in ordinary
geometry. This is extremely well suited to study the global aspects of
geometries which in some sense resemble compact spaces with
positive-definite metrics. To study noncommutative analogs of
noncompact manifolds with metrics of arbitrary signature it is perhaps
more practical to use calculi which are based on sets of
derivations. We shall use this method here. In all cases the entire
calculus can be considered as implicit in the module structure of the
set of 1-forms. We shall consider only the cases where this module is
free as a left or right module. It will in general not be free as a
bimodule.

There are basically two points of view. One can start with a set of
derivations in the strict sense of the word, a set of linear maps of
the algebra into itself which satisfy the Leibniz rule and use them as
basis for the construction of the associated differential forms. Or
one can start with a set of differential forms obtained for example
from some covariance criterion and construct a set of possibly twisted 
derivations which are dual to the forms. By `twisted' here we mean
derivations which satisfy a modified form of the Leibniz rule. We
shall describe both points of view and compare them.

Let $\c{A}$ be an algebra and $\lambda_a$, $1 \leq a \leq n$ a set of
$n$ elements of $\c{A}$ which is such that only the identity commutes
with it. This rule implies that only multiples of the identity will have
a vanishing differential. We have obviously therewith excluded
commutative algebras from consideration. In the example we consider this
condition will not be satisfied, which explains why we can have a
noncommutative geometry with only one dimension. We shall comment on
this latter. We introduce a set of derivations $e_a$ defined on an
arbitrary element $f \in \c{A}$ by $e_a f = [\lambda_a, f]$.  We have
here given the $\lambda_a$ the physical dimensions of mass; we set this
mass scale equal to one. Suppose that the algebra is generated formally
by $n$ elements $x^i$. If one defines the differential of $f \in \c{A}$
by $df(e_a) = e_a f$ exactly as one does in ordinary geometry, or by any
other method, then one finds that in general
$$ 
dx^i (e_a) \neq \delta^i_a.
$$
The `natural' basis $e_a$ of the derivations are almost never dual to
the `natural' basis $dx^i$ of the 1-forms.  There are basically two ways
to remedy the above default. One can try to construct a new basis
$\theta^a$ which is dual to the basis of the derivations or one can
introduce derivations $\partial_i$ which satisfy a modified form of
the Leibniz rule and which are dual to the $dx^i$. One has then either,
or both, of the following equations:
$$
\theta^a (e_b) = \delta^a_b, \qquad dx^i (\partial_j) = \delta^i_j.
$$
In general these two points of view are equivalent.  By construction
the $\theta^a$ commute with all elements of the algebra. These
commutation relations define the structure of the 1-forms as a
bimodule over the algebra.

We recall briefly the construction based on
derivations~\cite{DimMad96}. One finds that for the `frame' or
`Stehbein' $\theta^a$ to exist the $\lambda_a$ must satisfy a
constraint equation
\be
2 \lambda_c \lambda_d P^{cd}{}_{ab} - 
\lambda_c F^c{}_{ab} - K_{ab} = 0                             \label{fund}
\ee
with all the coefficients lying in the center $\c{A}$. The first set
of coefficients must be non-vanishing if the module of 2-forms is to be
nontrivial; it is related to a quantity which satisfies a sort of
Yang-Baxter equation.  Equation~(\ref{fund}) gives to the set of
$\lambda_a$ the form of a twisted Lie algebra with a central
extension. It is obviously a very severe restriction. If the algebra
is a $*$-algebra then the $\lambda_a$ must be antihermitian if the
derivations are to be real. The involution can be extended to the
general forms as well as to the tensor product of 1-forms by
introducing a set $J^{ab}{}_{cd}$ of central elements. If one
introduces a covariant derivative and requires that it be real then
the left and right Leibniz rules are connected through the
$J^{ab}{}_{cd}$. If the $J^{ab}{}_{cd}$ satisfy the Yang-Baxter equation
then the extension of the covariant derivative to the tensor product of
two 1-forms is real. More details of this can be found
elsewhere~\cite{FioMad98a}.

The dual point of view~\cite{WesZum90} consists in choosing the
differentials $dx^i$ as the starting point and constructing from them
a set of twisted derivations which satisfy a modified Leibniz
rule. Although at first sight this method seems to be less general
than the first, being normally restricted to quantum spaces invariant
under the coaction of some quantum group, in fact, as we saw above,
the quantum-group structure is more or less implicit also in the first
approach in the form of the Yang-Baxter equation. The dual point of
view has also the advantage in the fact that the twisted derivations
can be given a bimodule structure and an associated phase space is
perhaps more naturally constructed.

A metric on an algebra $\c{A}$ can be defined~\cite{DubMadMasMou96} in
terms of the 1-forms of a differential calculus $\Omega^*(\c{A})$
as a bilinear map
\be
g: \Omega^1(\c{A}) \otimes_{\c{A}} \Omega^1(\c{A}) 
\rightarrow \c{A}                                       \label{contrametric}
\ee
or~\cite{Cas95} in terms of the twisted derivations $\c{X}$ as a 
bilinear map
\be
g^\prime: \c{X} \otimes_{\c{A}} \c{X} \rightarrow \c{A}. \label{cometric}
\ee
We have distinguished here the two maps but in the case which interests
us here they are essentially one and the same.  In terms of the basis
these equations can be written respectively as
\be
g(\theta^a \otimes \theta^b) = g^{ab}, \qquad
g^\prime (\partial_i \otimes \partial_j) = 
g^\prime_{ij}                                       \label{metriccomp}
\ee
Since the $\theta^a$ commute with the elements of the algebra one sees
from the sequence of identities
\be
f g^{ab} = g(f \theta^a \otimes \theta^b) = 
g(\theta^a \otimes \theta^b f) = g^{ab} f             \label{componcom}
\ee
for arbitrary $f \in \b{R}_q^1$ that the $g^{ab}$ must lie in the
center of $\b{R}_q^1$; they must be real numbers.
Since the $\partial_i$ do not commute with the elements of the algebra
one sees from the sequence of (in)equalities
$$
f g^\prime_{ij} = g^\prime(f \partial_i \otimes \partial_j) \neq 
g^\prime(\partial_i \otimes \partial_j f) = g^\prime_{ij} f
$$
for arbitrary $f \in \b{R}_q^1$ that the $g^\prime_{ij}$ cannot lie in the
center of $\b{R}_q^1$. The commutation relations between $f$ and
$g^\prime_{ij}$ are however in principle calculable in terms of the
commutation relations between $f$ and $\partial_i$. A more
detailed exposition of the geometry of the algebra $\b{R}_q^3$
has been given elsewhere~\cite{FioMad98b}. 

Suppose that one particular `coordinate' $x^i$ has a discrete spectrum
$\ket{k}$. Then it is possible to give an observational definition of
the distance $ds(k)$ between $\ket{k}$ and $\ket{k+1}$ in terms of $g$
or $g^\prime$ by identifying $dx$ as the difference between the two
corresponding eigenvalues. It is our main purpose to study this
identification in detail in a particularly simple case.

\sect{The $q$-deformed euclidean spaces}

The $q$-deformed euclidean spaces~\cite{FadResTak89} $\b{C}_q^n$ and
$\b{R}_q^n$ are algebras which are covariant under the coaction of the
quantum groups $SO_q(n)$. To describe them it is convenient to
introduce the projector decomposition of the corresponding braid
matrix
$$
\hat R = q P_s - q^{-1} P_a + q^{1-n} P_t
$$
where the $P_s$, $P_a$, $P_t$ are $SO_q(n, \b{R})$-covariant
$q$-deformations of respectively the symmetric trace-free,
antisymmetric and trace projectors. They are mutually orthogonal and
their sum is equal to the identity:
$$
P_s + P_a + P_t  = 1.
$$
The trace projector is 1-dimensional and its matrix elements can be
written in the form
$$
P_t{}^{ij}{}_{kl} = (g^{mn}g_{mn})^{-1} g^{ij}g_{kl},
$$
where $g_{ij}$ is the $q$-deformed euclidean metric.
The $q$-euclidean space is the formal associative 
algebra $\b{C}_q^n$ with generators $x^i$ and relations
$$
P_a{}^{ij}{}_{kl} x^k x^l=0
$$
for all $i,j$. One obtains the real $q$-euclidean space by
choosing $q \in \b{R}^+$ and by giving the algebra an involution
defined by 
$$
x_i^* = x^j g_{ji}.
$$
This condition is an $SO_q(n, \b{R})$-covariant condition and $n$
linearly independent, real coordinates can be obtained as combinations
of the $x^i$. The `length' squared
$$
r^2 := g_{ij} x^i x^j = x_i^* x^i
$$ 
is $SO_q(n,\b{R})$-invariant, real and generates the center
$\c{Z}(\b{R}_q^n)$ of $\b{R}_q^n$. We can extend $\b{R}_q^n$
by adding to it the square root $r$ of $r^2$ and the inverse
$r^{-1}$. For reasons to become clear below when we introduce
differential calculi over $\b{R}_q^n$ we add also an extra generator
$\Lambda$ called the dilatator and its inverse $\Lambda^{-1}$ chosen
such that
\be
x^i \Lambda = q \Lambda x^i.                                  \label{lambdax}
\ee
We shall choose $\Lambda$ to be unitary. Since $r$ and $\Lambda$ do
not commute the center of the new extension is trivial.

We shall be here interested only in the case $n=1$.  The algebra
$\b{R}_q^1$ has only two generators $x$ and $\Lambda$ which
satisfy the commutation relation $x \Lambda = q \Lambda x$. We shall
choose $x$ hermitian and $q \in \b{R}^+$ with $q > 1$. This is a
modified version of the Weyl algebra with $q$ real instead of with
unit modulus. We can represent the algebra on a Hilbert space
$\c{R}_q$ with basis $\ket{k}$ by
\be
x \ket{k} = q^k \ket{k}, \qquad  \Lambda \ket{k} = \ket{k+1}.  \label{rep}
\ee
This is an infinite-dimensional version of the basis introduced by
Schwinger~\cite{Sch60} to study the Weyl algebra when $q$ is a root of
unity.  It explains the origin of the expression `dilatator'. Contrary
to the case considered by Schwinger however the spectrum of $\Lambda$
is continuous.

Introduce the element $y$ by the action
\be
y \ket{k} = k \ket{k}                                             \label{y}
\ee
on the basis elements. Then the commutation relations between
$\Lambda$ and $y$ can be written as
\be
\Lambda^{-1} y \Lambda = y + 1.                            \label{lambda-y}
\ee
We can write $x = q^y$ as an equality within $\b{R}^1_q$. 
We shall on occasion renormalize y. We introduce a renormalization 
parameter $z$ as 
$$
z = q^{-1}(q-1) > 0.
$$
The renormalization is then given by the substitution
\be
z y \mapsto y.                                               \label{renorm}
\ee
With the new value of $y$ the spacing between the spectral lines
vanishes with $z$. We shall refer to the old units as Planck units and
the new ones as laboratory units.  One can show~\cite{Ped79} that the
von~Neumann algebra generated by $\Lambda$ and $x$ or $y$ is a factor
of type $\mbox{I}_\infty$.

\sect{The $q$-deformed calculi}

One possible differential calculus over the algebra $\b{R}_q^1$ is
constructed by setting $d\Lambda = 0$ and
$$
x dx = dx x, \qquad dx \Lambda = q \Lambda dx.
$$
The frame is given by $\theta^1 = x^{-1} dx$. This calculus has an
involution given by $(dx)^* = dx^*$ but it is not based on derivations
and it has no covariance properties with respect to $SO_q(1)$.

We consider therefore another differential
calculus $\Omega^*(\b{R}_q^1)$ based on the relations~\cite{CarSchWatWei91}
\be
x dx = q dx x, \qquad dx \Lambda = q \Lambda dx                \label{calculus}
\ee
for the 1-forms. If we choose 
$$
\lambda_1 = - z^{-1} \Lambda
$$
then 
$$
e_1 x = q \Lambda x, \qquad e_1 \Lambda = 0
$$ 
and the calculus~(\ref{calculus}) is defined by the condition 
$df (e_1) = e_1 f$ for arbitrary $f \in \b{R}_q^1$. By setting 
$$
\lambda_2 = z^{-1} x
$$
and introducing a second derivation
$$
e_2 \Lambda = q \Lambda x, \qquad e_2 x = 0
$$ 
one could extend the calculus~(\ref{calculus}) by the condition 
$df (e_2) = e_2 f$ for arbitrary $f \in \b{R}_q^1$. One would find 
$x d\Lambda = q d\Lambda x$. We shall not do so since it will be seen
that $\Lambda$ is in a sense an element of the phase space associated
to $x$ and we are interested in position-space geometry.

The adjoint derivation $e_1^\dagger$ of $e_1$ is defined by
$$
e_1^\dagger f = (e_1 f^*)^*.
$$
The $e_1^\dagger$ on the left-hand side is not an adjoint of an
operator $e_1$. It is defined uniquely in terms of the involution of
$\b{R}^1_q$ whereas $e_1$ acts on this algebra as a vector space.

Since $\Lambda$ is unitary we have $(\lambda_1)^* \neq - \lambda_1$
and $e_1$ is not a real derivation. We introduce~\cite{CarSchWatWei91}
therefore a second differential calculus $\bar\Omega^*(\b{R}_q^1)$
defined by the relations
\be
x \bar dx = q^{-1} \bar dx x, \qquad 
\bar dx \Lambda = q \Lambda \bar dx                         \label{barcalculus}
\ee
and based on the derivation $\bar e_1$ formed using 
$\bar \lambda_1 = - \lambda_1^*$. This calculus is defined by the
condition $\bar df (\bar e_1) = \bar e_1 f$ for arbitrary 
$f \in \b{R}_q^1$. The derivation $\bar e_1$ is also not real. It is
easy to see however that
\be
e_1^\dagger = \bar e_1                                     \label{derstarbar}
\ee
and therefore that $(df)^* = \bar df^*$. By simple induction we find 
that for arbitrary integer $n$
$$
e_1 x^n = z^{-1} (q^n - 1) \Lambda x^n, \qquad
\bar e_1 x^n = z^{-1} (1 - q^{-n}) \Lambda^{-1} x^n.
$$

We can represent also $\Omega^*(\b{R}^1_q)$ and $\bar\Omega^*(\b{R}^1_q)$
on $\c{R}_q$.  For the two elements $dx$ and $\bar dx$ we have 
respectively
\be
dx \ket{k} = \alpha q^{k+1} \ket{k+1}, \qquad
\bar dx \ket{k} = \bar \alpha q^k \ket{k-1}                 \label{diffrep}
\ee
with two arbitrary complex parameters $\alpha$ and $\bar \alpha$. One
sees that $(dx)^* = \bar dx$ if and only if $\alpha^* = \bar \alpha$.
It is possible to represent $d$ and $\bar d$ as the operators
$$ 
d = - z^{-1} \alpha \, \ad \Lambda, \qquad
\bar d = z^{-1}\bar\alpha \, \ad \Lambda^{-1}.
$$
It is easy to see that the commutation relations~(\ref{calculus})
and~(\ref{barcalculus}) are satisfied. The above representations are
certainly not unique~\cite{Sch98}. 

The frame elements $\theta^1$ and $\bar \theta^1$ dual to the
derivations $e_1$ and $\bar e_1$ are given by
\be
\begin{array}{lll}
\theta^1 = \theta^1_1 dx, &&\theta^1_1 = \Lambda^{-1} x^{-1}, \\
\bar\theta^1 = \bar\theta^1_1 \bar dx, &&
\bar\theta^1_1 = q^{-1} \Lambda x^{-1}.
\end{array}                                                     \label{frame}
\ee
On $\c{R}_q$ they become the operators
\be
\theta^1 = \alpha , \qquad \bar \theta^1 = \bar \alpha            \label{int}  
\ee
proportional to the unit element. They were so constructed. 
The algebra $\b{R}^1_q$ is a subalgebra of the graded algebra of forms
$\Omega^*(\b{R}^1_q)$ and the representation~(\ref{rep}) can be extended
to a representation of the latter. In fact since $\Omega^1(\b{R}^1_q)$
and $\bar\Omega^1(\b{R}^1_q)$ are free $\b{R}^1_q$-modules of rank one 
with respectively the special basis $\theta^1$ and $\bar \theta^1$ we
can identify
$$
\Omega^*(\b{R}^1_q) =
\textstyle{\bigwedge^*} \otimes \b{R}^1_q, \qquad
\bar\Omega^*(\b{R}^1_q) =
\textstyle{\bigwedge^*} \otimes \b{R}^1_q
$$
where $\bigwedge^*$ is the exterior algebra over $\b{C}^1$ and so the
extension is trivial.

{From} the two differential calculi $\Omega^*(\b{R}_q^1)$ and $\bar
\Omega^*(\b{R}_q^1)$ we would like to construct a real differential
calculus $\Omega_R^*(\b{R}_q^1)$ with a differential $d_R$ such that
$(d_R f)^* = d_R f^*$. The construction has nothing to do with the
structure of $\b{R}^1_q$ so we give it in terms of a general algebra
$\c{A}$.

\sect{A Real calculus}

Consider an algebra $\c{A}$ with involution over which there
are two differential calculi $(\Omega^*(\c{A}), d)$ and
$(\bar\Omega^*(\c{A}), \bar d)$ neither of which is necessarily
real. Consider the product algebra $\tilde \c{A} = \c{A} \times \c{A}$
and over $\tilde \c{A}$ the differential calculus
\be
\tilde \Omega^*(\tilde \c{A}) = 
\Omega^*(\c{A}) \times \bar\Omega^*(\c{A}).               \label{prod-cal}
\ee
It has a natural differential given by $\tilde d = (d, \bar d)$. The 
embedding
$$
\c{A} \hookrightarrow \tilde \c{A}
$$
given by $f \mapsto (f, f)$ is well defined and compatible with the
involution 
\be
(f,g)^* = (g^*,f^*)                                       \label{alg-invol}
\ee
on $\tilde \c{A}$.

Let $X$ and $\bar X$ be two derivations of $\c{A}$. Then 
$\tilde X = (X, \bar X)$ is a derivation of $\tilde \c{A}$. We recall
that a derivation $X$ of an algebra $\c{A}$ is real if for arbitrary
$f \in \c{A}$ we have $X f^* = (X f)^*$. We saw in the previous
section that $e_1$ and $\bar e_1$ are not real. Then $\tilde X$ is a
real derivation if
\be
\tilde X (f,g))^* = (\tilde X (f, g))^*.                   \label{real-der}
\ee
This can be written as the conditions 
$$
\bar X f^* = (X f)^*, \qquad X g^* = (\bar X g)^*.
$$
The essential point to notice is that $\c{A}$ does not necessarily
remain invariant under real derivations of $\tilde \c{A}$. This is to
be expected since if $\c{A}$ had `interesting' real derivations they
could be used to construct directly a real differential calculus over
$\c{A}$.  

Suppose that $\Omega^*(\c{A})$ is defined in term of a set of inner
derivations $e_a = \ad \lambda_a$ and that $\bar\Omega^*(\c{A})$ is
defined in term of a set of inner derivations 
$\bar e_a = \ad \bar \lambda$. Suppose also that the corresponding 
$\tilde e_a = (e_a, \bar e_a)$ are real derivations of $\tilde \c{A}$.
{From} (\ref{real-der}) we see that this will be the case if and only if
$\bar \lambda_a = - \lambda_a^*$. We saw in the previous section that
$\tilde e_1 = (e_1, \bar e_1)$ is real and that in fact 
$\bar \lambda_1 = - \lambda_1^*$.  We define an involution on 
$\tilde \Omega^*(\tilde \c{A})$ by the condition
$$
(\tilde d (f, g))^* (\tilde e_a) = (\tilde e_a (f,g))^* = 
\tilde e_a (g^*,f^*).
$$
The differential $\tilde d$ is real by construction~\cite{FioMad98a}.

Define $\c{A}_R$ to be the smallest algebra which contains $\c{A}$ and
which is stable under the action of the derivations $\tilde e_a$. The
image in $\tilde \c{A}$ of the commutative subalgebra 
$\c{A}_0 \subset \c{A}$ of $\c{A}$ generated by $x$ is invariant under
the involution (\ref{alg-invol}). Define $e_{Ra}$ to be the
restriction of $\tilde e_a$ to $\c{A}_R$ and $d_R$ to be the
restriction of $\tilde d$ to $\c{A}_R$.  We have then
\be
d_R f (e_{Ra}) = (e_a f, \bar e_a f)                            \label{real-d}
\ee
and $d_R$ is also real.  We define 
\be
\Omega_R^1(\c{A}) \subset \tilde \Omega^1(\tilde \c{A})       \label{submodule}
\ee
to be the $\c{A}_R$-bimodule generated by the image of $d_R$. We write
$\Omega_R^1(\c{A})$ instead of $\Omega_R^1(\c{A}_R)$ since we keep 
$\Omega_R^0(\c{A}) = \c{A}$. The module structure determines a
differential calculus $(\Omega_R^*(\c{A}), d_R)$.  Suppose there
exists a frame $\theta^a$ for $\Omega^*(\c{A})$ and a frame 
$\bar \theta^a$ for $\bar \Omega^*(\c{A})$. We can extend also the
involution (\ref{alg-invol}) to all of $\Omega_R^*(\c{A})$ by setting
$$
(\theta^a)^* = \bar \theta^a
$$
and we can define $\Omega_R^1(\c{A})$ to be the $\c{A}_R$ module
generated by
$$
\theta_R^a = (\theta^a, \bar\theta^a).
$$
This is consistent with the previous definition since
$$
d_R f = e_{Ra} f \theta_R^a, \qquad e_{Ra} f \in \c{A}_R.
$$
{From} the relations
\be
\begin{array}{ll}
\theta^a(e_b) = \delta^a_b,      &\theta^a(\bar e_b) = 0, \\
\bar \theta^a(e_b) = 0, &\bar\theta^1(\bar e_b) = \delta^a_b    \label{array}
\end{array}
\ee
it follows that the frame dual to the derivation $e_{Ra}$ is indeed
$\theta_R^a$:
$$
\theta_R^a (e_{Rb}) = \delta^a_b.
$$
If we define the `Dirac operators' 
\be
\theta = - \lambda_a \theta^a, \qquad
\bar\theta = - \bar\lambda_a \bar\theta^a, \qquad
\theta_R = - \lambda_{Ra} \theta_R^a                              \label{Dirac}
\ee
then we find from the Equation~(\ref{real-d}) that for all $f \in \c{A}$
$$
d f = - [\theta, f], \qquad
\bar d \bar f = - [\bar\theta, f],\qquad
d_R f = - [\theta_R, f].
$$
Except for $\Omega_R^0(\c{A}) = \c{A}$ we can write
$$
\Omega^*_R(\c{A}) = \textstyle{\bigwedge}^* \otimes \c{A}_R
$$
where $\bigwedge^*$ is the algebra over $\b{C}$ generated by the
$\theta_R^a$.

We are now in a position to construct a real differential calculus
over $\b{R}^1_q$. According to the general remarks we see that 
$e_{R1} = (e_1, \bar e_1)$ is a real derivation of $\b{R}^1_{qR}$ and 
and that it is inner
\be
e_{R1} = \ad \lambda_{R1} \qquad
\lambda_{R1} = (\lambda_1, \bar\lambda_1) =
z^{-1} (- \Lambda, \Lambda^{-1}).                               \label{lambdaR}
\ee
Because of the identity
\be
e_{R1} x = (q \Lambda, \Lambda^{-1}) x                      \label{eRx}
\ee
we conclude that
\be
x d_R x = (q, q^{-1}) d_R x x, \qquad 
d_R x \Lambda = q \Lambda d_R x.                           \label{realcalculus}
\ee
These are the real-calculus equivalent of the relations~(\ref{calculus})
and~(\ref{barcalculus}). A representation of the 1-forms of the
differential calculus $\Omega_R^*(\b{R}^1_q)$ can be given on the
direct sum $\c{R}_q \oplus \c{R}_q$ of two separate and distinct
copies of $\c{R}_q$, one for $dx$ and one for $\bar dx$.
{From}~(\ref{diffrep}) one sees that $d_R x$ can be represented by the
operator
$$
d_R x \ket{k} = q^k (q \alpha \ket{k+1} 
+ \bar \alpha \overline{\ket{k-1}}).                     
$$
We have placed a bar over the second term to underline the fact that it
belongs to the second copy of $\c{R}_q$. 

Since the Equations~(\ref{array}) involve (in the case 
$\c{A} = \b{R}^1_q$) $e_1$ and $\bar e_1$ considered as derivations they
cannot be implemented on $\c{R}_q$.  However $e_1$ and $\bar e_1$ can be
considered as `annihilation' operators which map $\Omega_R^1(\b{R}^1_q)$
into $\Omega_R^0(\b{R}^1_q)$. Similarly $\theta^1$ and $\bar \theta^1$
have an interpretation~\cite{KasMadMas97} as `creation' operators which
take $\Omega_R^1(\b{R}^1_q)$ into $\Omega_R^2(\b{R}^1_q) \equiv 0$.  On
$\c{R}_q \oplus \c{R}_q$ the involution is given by the map 
$\alpha \mapsto \bar\alpha$. We shall choose
\be
\alpha = 1, \qquad \bar \alpha = 1                           \label{intvalue}
\ee
so that the map simply exchanges the two terms of $\c{R}_q \oplus \c{R}_q$.
On $\c{R}_q \oplus \c{R}_q$ we have the representation
\be 
\theta_R^1 = 1.                                               \label{realint}
\ee

If $d_R$ is to be a differential then the extension to higher order
forms much be such that $d_R^2 = 0$. We have then
\be
(d_R x)^2 = 0.                                                  \label{nil} 
\ee
It follows that
\be
d_R \theta_R^1 = 0, \qquad (\theta_R^1)^2 = 0.          \label{Omega_R-struc}
\ee
The module structure of $\Omega_R^1(\b{R}_q^1)$ is given by the relations
(\ref{realcalculus}), which are equivalent to the condition that
$\theta^1_R$ commute with all the elements of $\b{R}^1_q$. 
The algebraic structure of $\Omega_R^*(\b{R}_q^1)$ is defined by the
relations (\ref{Omega_R-struc}).

The algebra $\b{R}^1_q$ is a subalgebra of the graded algebra of forms
$\Omega^*(\b{R}^1_q)$ and the representation~(\ref{rep}) can be extended
to a representation of the latter. Again since $\Omega_R^1(\b{R}^1_q)$
is a free $\b{R}^1_q$-module of rank one with the special basis
$\theta_R^1$ we can identify
$$
\Omega_R^*(\b{R}^1_q) =
\textstyle{\bigwedge^*} \otimes \b{R}^1_q
$$
where $\bigwedge^*$ is the exterior algebra over $\b{C}^1$ and so the
extension is trivial.  The $\theta^1_R$ here is to be interpreted as an
element on the $\bigwedge^*$ and the equality gives its representation
as the unit in $\b{R}^1_q$. The second of
Equations~(\ref{Omega_R-struc}) is to be interpreted then as the
equation $1 \wedge 1 = 0$ in the exterior algebra.

The forms $\theta^1$, $\bar \theta^1$ and $\theta_R^1$ are
closed. They are also exact. In fact if we define 
$K \in \b{R}^1_q \times \b{R}^1_q$ by
\be
K = z (\Lambda^{-1}, \Lambda), \qquad K^* = K,                     \label{K}
\ee
then we find that
$$
\theta^1 = d (z\Lambda^{-1} y), \qquad
\bar \theta^1 = \bar d (z \Lambda y), \qquad
\theta_R^1 =  d_R (Ky).
$$

One can always write a number $x$ as the sum of a complex number $z$
and its complex conjugate $\bar z$. If some invariance property were
to forbid us from writing any formula involving $dx$ then we would
have to express it in terms of $dz$ and $d\bar z$.  What we have done
in this section is equivalent to just this. It is not even interesting
from the point of view of module structure; we have considered the
simple direct sum of two free modules and the submodule defined
in~(\ref{submodule}) is also free, with $\theta^1_R$ as a
generator. To a certain extent what we have done is similar in spirit
to the doubling of the rank of the module of 1-forms proposed by
previous authors accompanied by an `abstract' isomorphism~\cite{Fio94}
to then effectively reduce the rank by one half. One can also
construct a (smaller) real differential calculus over $\b{R}^1_q$ using
the derivation $\ad (\lambda_1 + \bar \lambda_1)$ but this calculus
has a set of 1-forms with a more complicated module structure than the
one we have constructed. There would be no simple relation like
(\ref{realcalculus}) between $xd_R x$ and $d_Rx x$. 

The limit $q\to 1$ is rather difficult to control. From the relations
of the algebra and the two differential calculi one might expect
$\Lambda \to 1$.  This is consistent with the limiting relations 
$e_1 x = \bar e_1 x = x$ and the intuitive idea that $x$ is an
exponential function on the line. However the representation 
(\ref{rep}) of the algebra becomes quite singular. In the
representation one has rather $x \to 1$. This would imply that the
parameters $\alpha$ and $\bar\alpha$ must tend to zero as $q\to 1$.
If one renormalizes according to (\ref{renorm}) then one finds that
the relation (\ref{lambda-y}) is consistent with the limit 
$\Lambda \to 1$ as $q \to 1$. We shall assume this to be the case.
We have then
$$
\lim_{q\to 1} \c{A}_R = \c{A}
$$
and the real differential calculus coincides with the diagonal
elements of the product in (\ref{prod-cal}).

\sect{The $q$-deformed derivatives}

We now look at the differential calculus from the dual point of
view. We introduce a twisted derivation $\px$ dual to the differential
$d$. For every $f \in \b{R}_q^1$ we require that $df(\px) = \px f$. If
one uses the (historical) convention of writing $df = dx f_1$, with
the differential to the left, then this means that for arbitrary 
$f \in \b{R}_q^1$
\be
\px f = dx(\px) f_1 = f_1.                                  \label{dualcond}
\ee
Consider the case $f = x^2$. Then $df = dx (1+q) x$ and so 
$f_1 = (1+q)x$. But $df(\px) = \px x^2$. Therefore $\px x^2 = (1+q) x$
Consider the case $f = \Lambda x$. Then $df = dx q^{-1} \Lambda$ and so 
$f_1 = q^{-1}\Lambda$. But $df(\px) = \px (\Lambda x)$. Therefore 
$\px (\Lambda x) = q^{-1} \Lambda$. By considering arbitrary polynomials in
$x$ and $\Lambda$ one finds the commutation relations
\be
\px x = 1 + q x \px, \qquad 
\px \Lambda = q^{-1} \Lambda \px.                     \label{comrel}
\ee
It is to be noticed here that the module structure of the differential
forms is considered as fixed and the commutation relations above are
derived from it. When a differential calculus is based on derivations
the module structure of the forms is derived from the Leibniz rule.
Notice also that we are here considering $x$ as an operator on
$\b{R}_q^1$ considered as a vector space. We should in principle put a
hat on it to distinguish it from the element $x$ in $\b{R}_q^1$
considered as an algebra. We are also considering $\px$ as an
operator on $\b{R}_q^1$ considered as a vector space. We should put a
hat on it also to distinguish it from the twisted derivation of
$\b{R}_q^1$.  We have effectively enlarged the algebra $\b{R}_q^1$ to
an algebra $T_q\b{R}_q^1$ by adding to it the element $\px$
with the commutation relations~(\ref{comrel}).

We noticed above that the differential $dx$ was not real. In general
$(df)^* \neq df^*$. Closely related to this is the fact that the
derivation $\px$ is not real. In general 
$(\px f)^* \neq \px f^*$. Therefore $\px$ is not 
antihermitian considered as an element of $T_q\b{R}_q^1$. 
One can introduce a second
twisted derivation $\bpx$ dual to the differential $\bar d$.
It is defined by the commutation relations
\be
\bpx x = 1 + q^{-1} x \bpx, \qquad 
\bpx \Lambda = q^{-1} \Lambda \bpx.                        \label{barcomrel}
\ee
We have then a second extension $\bar T_q\b{R}_q^1$ of $\b{R}_q^1$.

The representation~(\ref{rep}) of $\b{R}^1_q$ can be extended to a
representation of $T_q\b{R}_q^1$ and $\bar T_q\b{R}_q^1$. We have
respectively~\cite{SchWes92, HebSchSchWeiWes94}
\be
\begin{array}{l}
\px \ket{k} = - z^{-1} q^{-k-1} \ket{k} 
 + z^{-1}\beta q^{-k} \ket{k-1}, \\
\bpx \ket{k} = z^{-1} q^{-k} \ket{k} 
 + z^{-1}\bar\beta q^{-k-1} \ket{k+1}
\end{array}                                                   \label{derrep}
\ee
with again two arbitrary complex parameters $\beta$ and $\bar\beta$.  It
is easy to see that the commutation relations~(\ref{comrel})
and~(\ref{barcomrel}) are satisfied. Again the representations are
certainly not unique. We shall conclude from~(\ref{derrep}) that
$$
\bar T_q\b{R}_q^1 = T_q\b{R}_q^1.
$$

Due to the presence of the unit on the right-hand side of the
commutation relations~(\ref{comrel}) and~(\ref{barcomrel}) the relation
between $\px^*$ and $\bpx$ is not as simple as it was
in the case of the differentials. Adding the adjoint of
Equations~(\ref{comrel}) to Equations~(\ref{barcomrel}) yields the
commutation relations
$$
x (q \px^* + \bpx) = 
q (q \px^* + \bpx)x, \qquad
\Lambda (q \px^* + \bpx)  = 
q (q \px^* + \bpx) \Lambda.
$$
{From} this we can conclude that
\be
q \px^* + \bpx = c_1 \Lambda x^{-1}                        \label{starbar}
\ee
for some constant $c_1$. In terms of the parameters of the
representation~(\ref{derrep}) we find the expression
$$
c_1 = z^{-1}(\beta^* + q^{-1} \bar \beta) 
$$
for $c_1$.

We can consider the derivations $e_1$ and $\bar e_1$ also as
elements of $T_q\b{R}_q^1$. As such they satisfy the commutation 
relations
\be
\begin{array}{ll}
e_1 x = q \Lambda x + x e_1, 
&e_1 \Lambda = \Lambda e_1,\\
\bar e_1 x = \Lambda^{-1} x + x \bar e_1,   
&\bar e_1 \Lambda = \Lambda \bar e_1.                         \label{e-comrel}
\end{array}                                                   
\ee
These are the analogs of~(\ref{comrel}) and~(\ref{barcomrel}) respectively.
One sees immediately that as elements of $T_q\b{R}_q^1$ the $e_1$ and
$\bar e_1$ satisfy the relation
\be
e_1^* + \bar e_1 = c_2                                             \label{c2}
\ee
for some constant $c_2$. This is the derivation analog
of~(\ref{starbar}).  It is to be compared with~(\ref{derstarbar}).
The equation $e_1 f = [\lambda_1, f]$ relates the derivation $e_1$ to
the operator $\lambda_1$. There is an ambiguity
$$
\lambda_1 \mapsto \lambda_1 + z^{-1} \gamma
$$
in this identification which depends on a complex parameter $\gamma$. A
similar ambiguity exists for $\bar\lambda_1$. As
operators on $\c{R}_q$ we find that we can write then
\be
e_1 = z^{-1} \gamma + \lambda_1, \qquad
\bar e_1 = z^{-1} \bar\gamma + \bar \lambda_1                    \label{e1op}
\ee
and in terms of $\gamma$ and $\bar\gamma$ we find the expression
\be
c_2 = z^{-1} (\gamma^* + \bar \gamma)                         \label{c2gamma}
\ee
for $c_2$.

If we use the expressions~(\ref{frame}) then by
comparing~(\ref{comrel}) and~(\ref{barcomrel}) with~(\ref{e-comrel})
we deduce as above the relations
$$
\px - \Lambda^{-1} x^{-1} e_1 = c_3 \Lambda^{-1} x^{-1}, \qquad
\bpx - q^{-1} \Lambda x^{-1} \bar e_1 = c_4 \Lambda x^{-1},
$$
where $c_3$ and $c_4$ are two arbitrary constants.  We shall here set
these two constants equal to zero. This means that we choose
\be
\gamma = \beta, \qquad \bar\gamma = \bar\beta.           \label{gammabeta}
\ee
We find then the relations 
\be
\px = \theta^1_1 e_1 = \Lambda^{-1} x^{-1} e_1, \qquad
\bpx = \bar\theta^1_1 \bar e_1 = 
q^{-1} \Lambda x^{-1} \bar e_1                           \label{transfor}
\ee
between the derivations $(e_1, \bar e_1)$ and the twisted derivations
$(\px, \bpx)$. We recall that the vector space $\mbox{Der}(\b{R}^1_q)$
is not a left module over the algebra $\b{R}^1_q$. As operators on 
$\c{R}_q$ one finds the representations
$$
e_1 \ket{k} = - z^{-1} \ket{k+1} + z^{-1} \beta \ket{k}, \qquad
\bar e_1 \ket{k} = z^{-1} \ket{k-1} + z^{-1} \bar\beta \ket{k}
$$
for the derivations. It follows directly from~(\ref{e1op}) that
\be
\beta \Lambda^{-1} = 1 + q z x \px, \qquad
\bar\beta \Lambda = - 1 + z x \bpx.                   \label{lambda} 
\ee
Inverting these expressions, we find that in fact
$$
\bar T_q\b{R}_q^1 = T_q\b{R}_q^1 = \b{R}_q^1.
$$
Using the adjoint of the Equations~(\ref{comrel}) we can
write~(\ref{transfor}) in the form
$$
x \px^* - 1 = \Lambda e_1^*, \qquad
x \bpx = \Lambda \bar e_1,
$$
from which using~(\ref{c2}) we deduce that
\be
x (\px^* + \bpx) = 
1 + z^{-1} (\beta +\bar \beta) \Lambda.          \label{defect}
\ee
It is interesting to note that $d\px$ and $d\bpx$ are well-defined and
not equal to zero.

From~(\ref{eRx}) we are prompted to introduce the antihermitian
element $e_{R1}$ of $T_q\b{R}_q^1$ with the commutation relations
$$
[e_{R1}, x] = (q \Lambda, \Lambda^{-1}) x, \qquad 
[e_{R1}, \Lambda] = 0.
$$
{From} the definition of $e_{R1}$ as derivation one sees that the solution
is given by
\be
e_{R1} = \lambda_{R1} + c_R                                    \label{eR}
\ee
for some complex parameter $c_R$.  If $\beta = \pm 1$, 
$\bar \beta = \mp 1$ then
\be
c_R = \pm z^{-1} (1, -1)                                       \label{c22}
\ee
and using~(\ref{lambda}) one finds the relation
$$ 
e_{R1} = \pm (x\bpx, qx\px)
$$
between $e_{R1}$ on the one hand and $x\px$ and $x\bpx$ on the other.

It does not seem to be possible to construct a real metric on the
twisted derivations without ambiguity. The problem is
complicated by the fact that, whereas $\bar dx = (dx)^*$,
from~(\ref{defect}) one sees that $\bpx \neq - \px^*$.  It would be
natural to define, for example, $\pr$ by the condition
$$
d_R x(\pr) = 1.
$$
However it is easy to see that this is not possible since if $\pr$ is
to be antihermitian as an operator then $\pr x$ cannot be hermitian as
an element of the algebra and so cannot be set equal to one. One could 
make the choice~\cite{FicLorWes96}
$$
\pr = {1 \over 2} (\px - \px^*)
$$
or the choice 
$$
\bar\pr = {1 \over 2} (\bar\px^* - \bar\px)
$$
or any combination of the two. We find from the 
representation~(\ref{derrep}) that
$$
x \pr = {1\over 2} z^{-1} \beta (q^{-1} \Lambda^{-1} - \Lambda), \qquad
\pr x = {1\over 2} z^{-1} \beta (\Lambda^{-1} - q^{-1}\Lambda)
$$
from which we conclude that
$$
q \pr x = {1\over 2} (q+1) \beta \Lambda^{-1} + x \pr, \qquad 
\Lambda \pr = q \pr \Lambda.
$$
This is to be compared with~(\ref{comrel}) and~(\ref{barcomrel}).
In particular we find that as twisted derivation
$$
d_Rx(\pr) = \pr x = {1\over 2q} (q+1) \beta \Lambda^{-1} \neq 1.
$$
In view of this ambiguity we shall use the derivation $e_{R1}$ to define
hermitian differential operators.

\sect{Integration}

Because of~(\ref{int}), (\ref{intvalue}) and~(\ref{realint}) 
we define~\cite{Con94} the (definite) integrals to be the linear maps 
from $\Omega^1_q(\b{R}^1_q)$, $\bar\Omega^1_q(\b{R}^1_q)$ and 
$\Omega^1_R(\b{R}^1_q)$ into the complex numbers given by 
respectively
$$
\int f_1 \theta^1 = \tr(f_1), \qquad 
\int \bar f_1 \bar \theta^1 = \tr (\bar f_1), \qquad 
\int f_R \theta^1_R = \tr(f_R). 
$$
In the last expression the trace includes the sum of the components.  
Since the `space' is `noncompact' we have
$$
\int \theta^1 = \tr (1) = \infty.
$$
In all cases the integral of an exact form is equal to zero. For
example
$$
\int df = \int e_1 \theta^1 = \tr([\lambda_1, f]) = 0.
$$
This is in fact rather formal since it is possible for the commutator
of two unbounded operators to have a non-vanishing trace. The integral
$$
\int dx = \tr (e_1 x) = q \tr(\Lambda x) = 0
$$
but the integral
$$
\int d(\Lambda^{-1} x) = \int \Lambda^{-1} dx = q \tr (x) = 
\sum_{-\infty}^\infty q^k = \infty.
$$
We can interpret the trace as an inner-product on the algebra by
setting
$$
\langle f \vert g\rangle \equiv \int (f^* g) \theta^1 = \tr (f^* g).
$$
The trace defines a state which characterizes the representation we
are using. It follows immediately from the definition that an operator
which is hermitian as an element of the algebra is also hermitian with
respect to the inner-product.

\sect{The geometry}

It is now possible to give an intuitive interpretation of the
metric~(\ref{contrametric}) in terms of observables. One can think of
the algebra $\b{R}_q^1$ as describing a set of `lines' $x$ embedded in
a `plane' $(x, \Lambda)$ and defined by the condition 
`$\Lambda = \mbox{constant}$'. To within a normalization the unique metric
is given by
\be
g(\theta_R^1 \otimes \theta_R^1) = 1.                      \label{realmetric}
\ee
Using it we introduce the element 
\be
g^{\prime 11} = g(d_R x \otimes d_R x) = 
(e_{R1}x)^2 g(\theta_R^1 \otimes \theta_R^1) = (e_{R1}x)^2        \label{g11}
\ee
of the algebra. Then
\be
\sqrt{g^{\prime 11}} = e_{R1} x, \qquad
\Big(\sqrt{g^{\prime 11}}\Big)^* = \sqrt{g^{\prime 11}}.        \label{rootg}
\ee
We have a representation of $x$ and $d_Rx$ on the
Hilbert space $\c{R}_q$. In this representation the distance $s$ along
the `line' $x$ is given by the expression
\be
ds(k) = \Vert \sqrt{g^\prime_{11}} d_R x (\ket{k} + 
\overline{\ket{k}})\Vert                                     \label{lineelem}
\ee
with as usual $g^\prime_{11} = (g^{\prime 11})^{-1}$. This comes
directly from the original definition of $dx$ as an `infinitesimal
displacement'. Using~(\ref{realint}) we find that
$$
ds(k) = \Vert \, \ket{k} + \overline{\ket{k}}\, \Vert = 1.
$$
The `space' is discrete~\cite{Sny47} and the spacing between `points'
is uniform. The distance operator $s$ can be identified with the
element $y$ introduced in~(\ref{y}). This means that if we measure $y$
using laboratory units, introduced in Equation~(\ref{renorm}) then we
shall do the same with $s$. In these units then the distance between
neighboring `points' is given by
$$
ds(k) = z.
$$

If one forgets the reality condition then one can introduce the
hermitian metric $g$ with $g(\bar\theta^1 \otimes \theta^1) = 1$. One
finds then 
\be
g^{\prime 11} = g(\bar dx \otimes dx) = \bar e_1 x e_1 x = q^2 x^2 
                                                       \label{hermitianmetric}
\ee
and one concludes that 
\be
ds(i) = \Vert \sqrt{g^\prime_{11}} dx \ket{k}\Vert = q.        \label{dist}
\ee
One can also introduce the hermitian metric $g$ with $g(e_1 \otimes
e_1^*) = 1$. One finds then
$$
g^\prime_{11} = g(\px \otimes \px^*) = 
\Lambda^{-1} x^{-1} g(e_1 \otimes e_1^*) x^{-1} \Lambda = q^{-2} x^{-2}
$$
and one finds again the expression~(\ref{dist}) for the distance. 

If one neglects also hermiticity and introduces a metric $g$ with 
$g(\theta^1 \otimes \theta^1) = 1$ then one finds that
\be 
g^{\prime 11} = g(dx \otimes dx) = (e_1 x)^2 = (q\Lambda x)^2. \label{g11'}
\ee
Since we have defined a `tangent space' $T_q\b{R}^1_q$ and a
`cotangent space' $\Omega^1(\b{R}^1_q)$ it is of interest to interpret
the metric as a map
$$
\Omega_R^1(\b{R}_q^1) \buildrel g \over \longrightarrow T_q\b{R}^1_q.
$$
This corresponds to the `raising of indices' in ordinary geometry.
As such it can be defined as the map $g(\theta^1) = e_1$. A short
calculation yields that this is equivalent to
\be
g(dx) = g^{\prime 11} \px                                   \label{raising}
\ee
as it should be.  Although both $dx$ and $\px$ have been represented
on the same Hilbert space we cannot conclude that in this
representation the map~(\ref{raising}) is given by $g = 1$. That is,
as operators on $\c{R}_q$, we have $dx \neq g^{\prime 11} \px$. One
finds in fact that
$$
(dx - g^{\prime 11} \px) \ket{k} = q^2 z^{-1} q^{2k} \ket{k+2} +
(\alpha - \beta q z^{-1}) q^k \ket{k+1}.
$$

We define covariant derivatives $D$ and $\bar D$ on 
$\Omega^1(\b{R}^1_q)$ as maps 
$$
\Omega^1(\b{R}^1_q) \buildrel D \over \longrightarrow
\Omega^1(\b{R}^1_q) \otimes \Omega^1(\b{R}^1_q), \qquad
\bar\Omega^1(\b{R}^1_q) \buildrel \bar D \over \longrightarrow
\bar\Omega^1(\b{R}^1_q) \otimes \bar\Omega^1(\b{R}^1_q)
$$
which satisfy~\cite{DubMadMasMou96} left and right Leibniz rules.  The
metric-compatible, torsion-free connections are given by the covariant
derivatives
$$
D \theta^1 = 0, \qquad \bar D \bar\theta^1 = 0.
$$
These equations can be written also as
$$
D (dx) = q^2 \Lambda^2 x \theta^1 \otimes \theta^1, \qquad 
\bar D (\bar dx) = \Lambda^{-2} x \bar\theta^1 \otimes \bar\theta^1.
$$
The real torsion-free covariant derivative compatible with the real
metric is given by
\be 
D_R \theta_R^1 = 0.                                       \label{realcovder}
\ee
This can also be written in the form
$$
D_R(d_R x) = (q^2 \Lambda^2, \Lambda^{-2}) x 
\theta_R^1 \otimes \theta_R^1.
$$
The generalized flip $\sigma_R$ is given~\cite{DubMadMasMou95} by
$\sigma_R = 1$. This yields~\cite{FioMad98a} the involution
$$
(\theta_R^1 \otimes \theta_R^1)^* = \theta_R^1 \otimes \theta_R^1
$$
on the tensor product if the covariant derivative (\ref{realcovder})
is to be real:
$$
D_R\xi^* = (D_R\xi)^*.
$$
The geometry is `flat' in the sense that the curvature tensor defined
by $D_R$ vanishes. The interpretation is somewhat unsatisfactory however
here because of the existence of elements in the algebra which do not
lie in the center but which have nevertheless vanishing exterior
derivative. These elements play a relatively minor importance in the
geometry of the algebras $\b{R}^n_q$ for larger values of 
$n$~\cite{SchWatZum92}.

\sect{Yang-Mills fields}

Consider the algebra $\c{A}_q$ obtained by adding a time parameter 
$t \in \b{R}$ to $\b{R}^1_q$: $\c{A}_q = \c{C}(\b{R})\otimes\b{R}^1_q$.
The tensor product is understood to include a completion with respect
to the topologies.  Choose $\c{H}$ as the $\c{A}_q$-bimodule which is
free of rank $r$ as a left or right module and assume that is can
be considered as an $\b{R}^1_{qR}$-bimodule.  Introduce a differential
calculus over $\c{A}_q$ by choosing the ordinary de~Rham differential
calculus over the time parameter and $\Omega_R^*(\b{R}^1_q)$ over the
factor $\b{R}^1_q$. One defines a covariant derivative of 
$\psi \in \c{H}$ as a map
$$
\c{H} \buildrel D \over \longrightarrow
\Omega^1(\c{A}_q) \otimes \c{H}
$$
which satisfies the left Leibniz rule
$$
D(f\psi) = df \otimes \psi + f D\psi.
$$
We shall henceforth drop the tensor product symbol and write
$$
D\psi = dt D_t \psi + D_R\psi.
$$
We define 
$$
D_t \psi = (\pt + A_t) \psi.
$$
Since $\c{A}_q$ is an algebra with involution we can choose as gauge
group the set $\c{U}_q(r)$ of unitary elements of $M_r(\c{A}_q)$. A
gauge transformation $g \in \c{U}_q$ is a map
$$
\psi \mapsto \psi g, \qquad A \mapsto g^{-1} A g + g^{-1} dg
$$ 
which is independent of $\Lambda$. It is easy to
see~\cite{DubKerMad89, Mad95} that $\phi_R = A_R - \theta_R$
transforms under the adjoint action of the gauge group. We define then
\be
D_R \psi = - \theta_R \psi - \psi \phi_R.                      \label{der-psi}
\ee
This covariant derivative is covariant under the right action of the
gauge group and satisfies a left Leibniz rule. The covariant
derivative and the field strength transform as usual
$$
D \psi \mapsto (D \psi) g, \qquad
F \mapsto g^{-1} F g
$$
One can also write $D_R\psi = \theta^1 D_{R1} \psi$ with 
$D_{R1} = e_{R1} + A_{R1}$. The field strength can be written then
$$
F\psi = D^2\psi =  dt \, \theta_R^1 \psi F_{t1}
$$
with
$$
F_{t1} = \pt A_{R1} - e_{R1} A_t.
$$
When the gauge potential vanishes one has from (\ref{der-psi}) 
$$
D_R \psi = \theta_R^1 e_{R1} \psi.
$$

To form invariants we introduce the metric~(\ref{realmetric}).  We 
define the matter action $S_M$ by analogy with the commutative case:
\be
S_M = \tr \int dt (D_t \psi^* D_t \psi + 
D_{R1} \psi^* D_{R1} \psi).                                        \label{S_M}
\ee
The trace is here over the Lie algebra of the gauge group and over the
representation of the algebra $\b{R}^1_q$. We define also as usual the
Yang-Mills action $S_{YM}$ as
$$
S_{YM} = {1\over 4} \tr (F_{t1} F_{t1})
$$
and the action to be the sum $S = S_M + S_{YM}$. The trace however
would depend on the representation of the algebra and it is not
obvious how one should vary $S$.  To define the trace we must consider
explicitly the representations of $\psi$ and $A_t$ and
$A_{R1}$ on the Hilbert space $\c{R}_q$.  Since `space' has only one
dimension there are no dynamical solutions to the vacuum Yank-Mills
equations. There is no dispersion relation since there are no
transverse modes.  One can also write the action as an integral using
the definition of Section~6.

In the spirit of noncommutative geometry the `state vectors' play the
role of the set of points.  The eigenvalues of an observable of the
algebra, in a given representation, are the noncommutative equivalents
of the values which its classical counterpart can take. An eigenvector
associated to a given eigenvalue describes a set of states in which
the given observable can take the prescribed value. This is exactly
like quantum mechanics but in position space. Consider now a field
configuration, for example an element of the initial algebra $\c{A}$
if it is a scalar field or an element of an algebra of forms over
$\c{A}$ if it is a Yang-Mills field. Suppose that both of these
algebras have a representation on some Hilbert space and suppose that
there exists a well-defined energy functional which is also
represented as an operator on the Hilbert space.  A vacuum
configuration would be then an element of the algebra which is such
that the expectation value of the corresponding value of the energy
functional in any state vanishes. This is the same as saying that a
field is equal to zero if the value of its energy is equal to zero at
every point of space. 

In the `classical' noncommutative case a derivation, if it exists at
all, is a map of the algebra into itself; it is not an element of the
algebra. In the case we are considering this is not so. The algebra
$\b{R}^1_q$ is a position space described by the subalgebra generated
by $x$ extended by $\Lambda$ which is an element of the associated
phase space. The differential calculus however is somehow restricted
to the position space by the condition $d\Lambda = 0$.  Both the
initial algebra and the algebra of forms contain then operators which
correspond to derivations. We have in fact given the representation of
these elements on $\c{R}_q$, the same Hilbert space on which the
`position' variables and the forms are represented. A vacuum
configuration is then something different than it is in the
`classical' case. 

Consider, for example, a scalar field $\psi(x) \in \b{R}^1_q$ and
suppose that the energy functional is of the simple form 
$\c{E} = (e_{R1} \psi)^* (e_{R1} \psi)$. If $e_{R1}$ is considered as
partial derivative then $\c{E} = \c{E}(x)$ depends on the position
variable alone and a vacuum configuration would be one in which the
expectation value of $\c{E}$ vanishes for all state vectors. This would
normally be one with $\psi = \psi_0$ for some $\psi_0$ ith $e_{R1}
\psi_0 = 0$.  However $e_{R1}$ as operator belongs also to $\b{R}^1_q$
and the expression for the energy functional could be interpreted as one
quadratic in this element. In this case the only possible vacuum
configuration would be $\psi = 0$. There exist particular state vectors
for which the energy functional of more complicated configurations
vanish. As an example of this we return to the Yang-Mills case.  One
would like a vacuum to be given as usual by $\psi = 1$ (the unit cyclic
vector of $\c{H}$) and $A_R = 0$.  One finds then as condition that
$$
D_R\psi (\ket{k} + \overline{\ket{k}}) = 
d_R (\ket{k} + \overline{\ket{k}}) = 0.
$$
To be concrete we shall suppose that $c_R$ is given by~(\ref{c22}).
From~(\ref{lambdaR}) one sees that the vacuum equation leads to the
conditions
$$
(\Lambda - 1) \sum_k a_k \ket{k} = 0, \qquad
(\Lambda^{-1} - 1) \sum_k \bar a_k \ket{k} = 0
$$
on the two copies of $\c{R}_q$. The vacuum state vectors form then a
subspace of $\c{R}_q$ of dimension 2 spanned by the vectors given by
$a_k = 1$, $\bar a_k = 1$. These values depend of course on our
choice of $c_R$.  All vacuum state vectors have infinite norm. The
vacuum state vectors would be the analog of the vacuum of quantum
field theory which is defined as the vector in Fock space which is
annihilated by the energy-momentum vector of Minkowski space. The
Fock-space vector is taken to be of unit norm.

\sect{The Schr\"odinger equation}

Recall that on a curved manifold with metric $g_{\mu\nu}$ the laplacian
is defined to be the hermitian operator
$$
\Delta = - g^{\mu\nu}D_\mu D_\nu = 
- - {1 \over \sqrt g} \partial_\mu (\sqrt g g^{\mu\nu} \partial_\nu).
$$
Because of~(\ref{der-psi}) on the geometry defined by 
$\Omega_R^*(\b{R}^1_q)$, with metric~(\ref{realmetric}), this becomes
\be
\Delta_R \psi = - e_{R1}^2 \psi.                               \label{Laplace} 
\ee
We shall suppose that the gauge-covariant Schr\"odinger equation has
the usual form
\be
i D_t \psi = {1 \over 2m} \Delta_R \psi               \label{Schroedinger}
\ee
where $\Delta_R$ is the Laplace operator~(\ref{Laplace}).  There is a
conserved current which we write in the form
\be 
\pt \rho = D_{R1} J_R^1                                    \label{conserv}
\ee
with as usual
$$
\rho = \psi^* \psi, \qquad 
J_{R1} =  {i \over 2m} (\psi^* D_{R1} \psi - D_{R1} \psi^* \psi).
$$
The conservation law follows directly from the field
equations. Normally one derives the latter from an action
principle. In the present situation this would be a non-relativistic
form of the expression (\ref{S_M}):
$$
S = \tr \int dt (i \psi^* D_t \psi -
{1 \over 2m} D_{R1} \psi^* D_{R1} \psi).
$$

Consider the relativistic case and assume the usual form
\be
- - \pt^2 \psi = \Delta_R \psi + m^2\psi                         \label{K-G}
\ee
for the Klein-Gordan equation. Suppose that $A_R = 0$.  The laplacian
has then a set of `almost' eigenvectors. From the commutation
relations
$$
e^{iky} \Lambda = e^{ik} \Lambda e^{iky}
$$
one finds that
$$
e_1 e^{iky} = z^{-1} (e^{ik} - 1) \Lambda e^{iky}, \qquad
\bar e_1 e^{iky} = 
z^{-1} (1 - e^{-ik}) \Lambda^{-1} e^{iky}
$$
from which it follows that
$$
e_{R1} e^{iky} = ikL e^{iky}
$$
where 
$$
L = {1 \over 2 ikz} \Big((e^{ik}- 1)\Lambda,\;
(1 - e^{-ik})\Lambda^{-1}\Big).
$$
{From} the expression (\ref{Laplace}) one concludes then that
$$
\Delta_R e^{iky} = k^2 L^2 e^{iky}.
$$
We could renormalize the space unit as in Equation~(\ref{renorm}) to
laboratory units. If we keep the Planck units we must renormalize the
time unit so it will be also in Planck units. We do this by the
transformation
$$
z^{-1}t \mapsto t.
$$
We find then that 
$$
\psi = e^{- i\omega L t - ky)}
$$
is a solution to~(\ref{K-G}) provided the dispersion relation
$$
(\omega^2 - k^2) z^2 L^2 = m^2  
$$
is satisfied.  However the above dispersion relations are misleading
since $\omega$ can not be identified with the energy; the coefficient
of the time coordinate is in fact the product $\omega L$ and we must
set therefore
$$
E = \omega L.
$$

We would like to consider $\psi$ as an element of a free
$\c{A}_0$-module. We recall that $\c{A}_0$ is the commutative
subalgebra of $\c{A}$ generated by $x$. In general however $\c{A}_0$
is not invariant under the action of the hermitian derivations. We
consider then the limit $q\to 1$. In this limit we have argued that 
$\Lambda \to 1$ but at the same time $z \to 0$ so the following 
argument is subject to caution. We supposed that as $q \to 1$ we have
$\bar e_1 \to e_1$. In this rather singular limit we can identify then
$$
e_{R1} = \frac 12 (e_1 +\bar e_1)(1,1) + o(z)
$$
and in this limit
$$
L = z^{-1}{\sin k \over k} (1,1) + o(1).
$$
This second equality seems to follow from the first but it is
especially difficult to justify. If we accept it however then with the
new time unit we find that
$$
E^2 = \omega^2 {\sin^2 k \over k^2}
$$
and the dispersion relation in the relativistic case becomes
\be
E^2 = m^2 + \sin^2 k.                                         \label{E}
\ee
If $k = \pi n$, with $n \in \b{Z}$ then $E = 0$ and one has
$$
e_{R1} e^{-iky} = 0.
$$
In the massless case this yields a set of `stationary-wave' solutions
to the field equations.

When $k << \pi/2$ (in Planck units) one obtains the usual dispersion
relation $E^2 = m^2 + k^2$. In the case $m << 1$ as $k\to \pi/2$ then
$E$ tends to a maximum value equal to $1$ (again in Planck
units). Values of $k$ greater than $\pi/2$ would be difficult to
interpret physically.. For comparison we recall that, neglecting the
gap corrections, the dispersion relation for acoustical phonons on a
lattice~\cite{Kit96} is of the form
$$
E^2 = \sin^2 {k\over 2}.
$$
Here $E$ is the phonon energy and $k$ is the wave number.  This has
the same form as~(\ref{E}) when $m=0$. The factor $1/2$ is a
convention. The first Brillouin zone is the range $-\pi \leq k \leq
\pi$. There are also optical phonons which are similar to the case
$m>0$ but they have a different dispersion relation. The `space'
$\b{R}^1_q$ is not an ordinary crystal.

\sect{Phase space}

If we wish to construct a real phase `space' associated to the
position `space' we must define two hermitian operators which can play
the role of `position' and `momentum'.  We have already remarked that
the distance operator $s$ can be identified with the element $y$
introduced in~(\ref{y}). As `position' operator we choose then the
renormalized $y$ given by~(\ref{renorm}). A short calculation shows
that
$$
e_{R1} y =  z K^{-1}
$$
where the element $K$ was introduced in~(\ref{K}).  If we consider
then $e_{R1}$ as an operator we have the commutation relation
$$
[e_{R1}, y] = (\Lambda, \Lambda^{-1}).
$$
We have already noticed that $e_{R1}$ is antihermitian. We define then
the momentum associated to $y$ to be
$$
p_y = - i e_{R1} = i z^{-1} (\Lambda, - \Lambda^{-1}).
$$
We have not written the extra constant term $c_R$ of
Equation~(\ref{c22}) since it does not contribute to the commutation
relation:
\be
[p_y, y ] = - i h.                                              \label{cr}
\ee
We have here introduced the hermitian element
$$
h  = (\Lambda,  \Lambda^{-1})
$$
of $\b{R}^1_q \times \b{R}^1_q$. Since we suppose that $\Lambda \to 1$
as $q \to 1$ we see that the commutation relation (\ref{cr}) becomes the
ordinary one in this limit.

We introduce the `annihilation operator' 
\be
a = {1 \over \sqrt 2} (y + i p_y).                               \label{a}
\ee
Then from~(\ref{cr}) follows the commutation relation
\be
[a, a^*] = h.                                              \label{a-a^*}
\ee
It is not possible to express $h$ in terms of $a$ and $a^*$. The
operator $e_{R1}$ was taken as the antihermitian part of $e_1$; the
operator $h$ depends also on the hermitian part. From~(\ref{lambda-y})
we find however that
\be
[a,h] = {1 \over 2} z^2 (a^* - a).                         \label{a-h}
\ee
To define a vacuum and a number operator we must `dress' the operator
$a$, introduce an operator $b$ so that the standard relations
$[b,b^*]=1$ hold. It does not seem to be possible to do this exactly
but it can be done as a perturbation series in $z$. One finds from
(\ref{a-a^*}) and (\ref{a-h}) that
$$
b = h^{-1/2} a + \frac 14 z^2 a + \frac 16 z^2 (a - a^*)^3 + o(z^4).
$$
The vacuum is chosen then as usual by the condition $b\ket{0} = 0$,
the number operator is given by $N = b^* b$ and the number
representation $\ket{n}$ for $n \in \b{N}$ by
$$
\ket{n} = {1 \over \sqrt{n!}} (b^*)^n \ket{0}.
$$
{From} (\ref{a-h}) we find that
$$
[N,h] = \frac 12 z^2 ((b^*)^2 - b^2) + o (z^4).
$$

There have been several $q$-deformed versions of the harmonic
oscillator~\cite{Mac89, Bie89, CarSchWat91, Fio94, LorRufWes97}. We
shall reproduce here the equivalent calculations for the geometry
described in Section~5. As hamiltonian we choose
$$
H = {1\over 2} (\Delta_R + y^2)
$$
in Planck units. This can be written also as $H = a^* a + \frac 12 h$
and in terms of $b$ it is given by
$$
H = b^* h b + \frac 12 h - \frac 12 z^2 b^* b + 
\frac 16 z^2 ((b - b^*)^3 b - b^* (b - b^*)^3) + o(z^4).
$$
We see then that in terms of the `dressed' annihilation and creation
operators the `bare' hamiltonian is rather complicated. In particular the
`physical vacuum' is no longer an eigenvector of the `bare'
hamiltonian:
$$
H \ket{0} = \frac 12 \ket{0} + \frac 16 z^2 \ket{1} +  
{1 \over \sqrt 2} z^2 \ket{2} - 
{1 \over \sqrt 6} z^2 \ket{3} + o(z^4).
$$

\sect{Non-local metrics}

We have devoted special attention to one particular metric on the
calculus $\Omega^*_R(\b{R}^1_q)$ for reasons given in Section~1: it is
the only local metric. To test what this means in practice it is of
interest to examine other metrics, which necessarily do not fulfill
the locality condition. We would like the metric to have an associated
linear connection so we shall first examine the most general form
which this can take. We set as usual
$$
D_1\theta^1 = - \omega^1{}_{11} \theta^1 \otimes \theta^1, \qquad
\bar D_1 \bar\theta^1 = 
- - \bar\omega^1{}_{11} \bar\theta^1 \otimes \bar\theta^1
$$
as Ansatz for the linear connection. From the general 
theory~\cite{DubMadMasMou96} these must satisfy a left and right
Leibniz rule
$$
\begin{array}{ll}
D_1(f\theta^1) = df \otimes \theta^1 - 
f \omega^1{}_{11} \theta^1 \otimes \theta^1,
&D_1(\theta^1 f) = \sigma(\theta^1 \otimes df) -
\omega^1{}_{11} f \theta^1 \otimes \theta^1, \\
\bar D_1(f\bar\theta^1) = \bar df \otimes \bar\theta^1 - 
f \bar\omega^1{}_{11} \bar\theta^1 \otimes \bar\theta^1,
&\bar D_1(\bar\theta^1 f) = \bar\sigma(\bar\theta^1 \otimes \bar df) -
\bar\omega^1{}_{11} f \bar\theta^1 \otimes \bar\theta^1,
\end{array}
$$
where $f \in \b{R}^1_q$ and the generalized flips~\cite{DubMadMasMou95}
$\sigma$ and $\bar\sigma$ can be written as
$$
\sigma(\theta^1 \otimes \theta^1) = S \theta^1 \otimes \theta^1, \qquad
\bar\sigma(\bar\theta^1 \otimes \bar\theta^1) = 
\bar S \bar\theta^1 \otimes \bar\theta^1.
$$
{From} the compatibility conditions
$$
\begin{array}{ll}
D_1(\Lambda \theta^1) = D_1(\theta^1 \Lambda),
&D_1(x \theta^1) = D_1(\theta^1 x),             \\
\bar D_1(\Lambda \bar \theta^1) = \bar  D_1(\bar  \theta^1 \Lambda), 
&\bar D_1(x \bar\theta^1) = \bar D_1(\bar\theta^1 x)
\end{array}
$$
it is easy to see that, to within a multiplicative constant, there are
only two solutions, the one given previously in Section~6 and a new one
given by
\be
\omega^1{}_{11} = \Lambda,            \quad       S = q^{-1}, \qquad
\bar\omega^1{}_{11} = q \Lambda^{-1}, \quad  \bar S = q.  \label{lin-conn}
\ee
We set 
$$
g(\theta^1 \otimes \theta^1) = g^{11}, \qquad
g(\bar\theta^1 \otimes \bar\theta^1) = \bar g^{11}.
$$
The metric compatibility condition~\cite{DubMadMasMou95} can be written
$$
dg^{11} = - (1 + S) \omega^1{}_{11} g^{11} \theta^1, \qquad
\bar d \bar g^{11} = 
- - (1 + \bar S) \bar\omega^1{}_{11} \bar g^{11} \bar\theta^1.
$$
The possible solution to this equation, corresponding to the linear
connection~(\ref{lin-conn}), is given by
$$
g_{11} = (q\Lambda x)^2, \qquad \bar g_{11} = (\Lambda^{-1} x)^2.
$$
This can also be written in the form
\be
g(dx \otimes dx) = 1, \qquad 
g(\bar dx \otimes \bar dx) = 1                         \label{nonlocal}
\ee
and the corresponding covariant derivative can be written also as
$$
D_1 dx = 0, \qquad \bar D_1 \bar dx = 0.
$$
The `space' now is a discrete subset of the positive real axis with an
accumulation point at the origin.  The `non-locality' means that if
$f$ is a `function' and $\alpha$ a form then the norm of $f\alpha$
cannot be equal to $f$ times the norm of $\alpha$.  To see this we
multiple (\ref{nonlocal}) from the right by $x$. If we supposed that
the metric were left and right linear then we would find
$$
x = x g(dx \otimes dx) = g(x dx \otimes dx) = 
q^2 g(dx \otimes dx x) = q^2 g(dx \otimes dx) x = q^2 x.
$$
The first and fifth equalities are mathematical trivialities; the
third follows directly from (\ref{calculus}) . Therefore either the
second or forth, or both, must be false.  There are no 2-forms and so 
the curvature and torsion of the non-local metric vanish.

The Ansatz for a covariant derivative on the real calculus is
\be
D_{R1} \theta_R^1 = 
- - \omega^1_{R11} \theta_R^1 \otimes \theta_R^1.           \label{real-non-loc}
\ee
If the generalized flip is given by $\sigma_R = (q^{-1}, q)$ then
the appropriate involution~\cite{FioMad98a} on the tensor product is
given by 
$$
(\theta_R^1 \otimes \theta_R^1)^* = 
(q^{-1}, q) (\theta_R^1 \otimes \theta_R^1).
$$
The solution to (\ref{real-non-loc}) is given by 
$$
\omega^1_{R11} = (\Lambda, q \Lambda^{-1}), \qquad
(\omega^1_{R11})^* = (q, q^{-1})\omega^1_{R11}, 
$$
The connection coefficient is not hermitian but the covariant derivative
is real.

As an example of an application we return to the Yang-Mills fields
written, for simplicity, using the derivations $D_1$ instead of $D_R$ and
where now of course the $D_1$ must be chosen compatible with the new
metric~(\ref{nonlocal}).  The differential calculus over $\c{A}_q$ is
now the ordinary de~Rham differential calculus over the time parameter
and $\Omega^*(\b{R}^1_q)$ over the factor $\b{R}^1_q$. Otherwise all is
as before in Section~8 except that $A$ is now given by 
$$
A = dt A_t + \theta^1 A_1
$$
and 
$$ 
d\psi = dt \pt \psi + \theta^1 e_1 \psi
$$
with $e_1\psi = [\lambda_1, \psi]$. More important, the
action $S$ becomes
$$
S_M = \tr (D_t \psi^* D_t \psi) +  g^{11} \tr (D_1 \psi^* D_1 \psi)
+ {1\over 4} g^{11} \tr (F_{t1} F_{t1}).
$$
The metric coefficient $g^{11} = (q\Lambda x)^{-2}$ does not commute
with the other factors in this expression so there is an ordering
ambiguity. But it must be outside the trace in order not to destroy
gauge invariance. Motivated by Equation~(\ref{transfor}) we introduce
the `twisted' covariant derivative $\nabla_1$ by the 
equation
$$
D_1 = \sqrt{g_{11}}\, \nabla_1.
$$
If we write also $A_1 = \sqrt{g_{11}} \, B_1$ and set
$D_t = \nabla_t$ and $G_{t1} = \pt B_1 - \px B_t$ then we find that
$$
\nabla = \px + B_1, \qquad F_{t1} = \sqrt{g_{11}} \, G_{t1}
$$
and we can write the action in the form
$$
S_M = \tr (\nabla_t \psi^* \nabla_t \psi) +  
g^{11} \tr (\sqrt{g_{11}}\,\nabla_1 \psi^* \sqrt{g_{11}} \,\nabla_1 \psi)
+ {1\over 4} g^{11} \tr (\sqrt{g_{11}}\,G_{t1} \sqrt{g_{11}}\,G_{t1}).
$$
We identify $\psi$ with an element $\psi(x)$ in the subalgebra of
$\b{R}^1_q$ generated by the element $x$ and we suppose also that 
$B_1 = B_1(x)$. The action can be written then in the form
$$
S_M = \tr (\nabla_t \psi^* \nabla_t \psi) +  
\Lambda^{-1} \tr (\nabla_1 \psi^* \Lambda \nabla_1 \psi)
+ {1\over 4} \Lambda^{-1} \tr (G_{t1} \Lambda G_{t1}).
$$
It is the quantity $\Lambda G_{t1}$ which is gauge covariant.

As a second example of we return to the Schr\"odinger equation, written
again using the covariant derivative $D_1$ compatible with the
metric~(\ref{nonlocal}).  There are two possible forms for the Laplace
operator $\Delta$. In the absence of a gauge potential one can choose
either
$$
\Delta = - g^{11} D_1 D_1 = 
- - q \Lambda^{-2} x^{-2} e_1^2 + q \Lambda^{-1} x^{-2} e_1
$$
or
$$
\Delta = - {1 \over \sqrt g} (e_1 \sqrt g g^{11} e_1) =
- - q \Lambda^{-2} x^{-2} e_1^2 + \Lambda^{-1} x^{-2} e_1.
$$
The two coincide when $q=1$. We shall choose the latter. If we
introduce then the current `density'
$$
\sqrt g J^1 =
{i \over 2m} \Lambda^{-1} (\psi^* \Lambda \px  \psi - 
\Lambda \px  \psi^* \psi) 
$$
the right-hand side of~(\ref{conserv}) becomes
$$
{1 \over \sqrt g} e_1 (\sqrt g J^1) = 
{i \over 2m} q \Lambda^{-1} \px(\psi^* \Lambda \px \psi - 
\Lambda \px \psi^* \psi).
$$
We have here used the relations~(\ref{transfor}). The conservation law
becomes then 
$$ 
\pt \rho = {1 \over \sqrt g} e_1 (\sqrt g J^1).
$$
This is the equivalent of~(\ref{conserv}) in the new metric.

It is of interest to compare the structure of the `space' endowed with
the two different metrics we have considered. We saw that the weak
completion of the algebra $\b{R}^1_q$ was in both cases a
type-$\mbox{I}_\infty$ factor. The metric can have no effect on this
since the set of `points' is discrete and the induced measures are
absolutely continuous one with respect to the other.  With the first
metric the `space' is an equally spaced lattice structure within the
entire real line. With the second metric one finds a lattice structure
with variable spacing within the half-line $(0, \infty)$. In this case
it would be natural either to add the origin to obtain a `space' with
boundary or to add the origin and another copy of the `space' to
obtain again the entire real line. In either case the algebra is no
longer a factor. It is also to be noticed that the fact we obtained a
factor of type $\mbox{I}_\infty$ is due to a choice of representation
and not the structure of the algebra. Had we chosen a representation
with a continuous spectrum for $\Lambda$ the resulting factor would be
of type $\mbox{II}_\infty$~\cite{Sch98}.

\section*{Acknowledgments} The authors would like to thank G. Fiore, 
H. Grosse, L. Vainerman and K. Schm\"udgen for enlightening
discussions.  Two of them (RH) and (JM) would also like to thank the
Max-Planck-Institut f\"ur Physik in M\"unchen for financial support.

\end{document}